\documentclass[10pt]{amsart}

\newcommand{\BC}{{\mathbf C}}
\newcommand{\BF}{{\mathbf F}}
\newcommand{\BQ}{{\mathbf Q}}
\newcommand{\BQb}{\overline{\mathbf Q}}
\newcommand{\BR}{{\mathbf R}}
\newcommand{\BZ}{{\mathbf Z}}
\newcommand{\ba}{{\mathbf a}}
\newcommand{\bb}{{\mathbf b}}
\newcommand{\be}{{\mathbf e}}
\newcommand{\bx}{{\mathbf x}}

\newcommand{\pid}{{\mathfrak p}}
\newcommand{\Pid}{{\mathfrak P}}

\newcommand{\CI}{{\mathcal I}}
\newcommand{\CN}{{\mathcal N}}
\newcommand{\CO}{{\mathcal O}}

\newcommand{\eps}{\epsilon}
\newcommand\ra{\rightarrow}
\newcommand{\kb}{\overline{k}}
\newcommand{\Kp}{K^{+}}
\def\simp{{\rm simple}}
\def\asimp{{\rm abs.simple}}

\DeclareMathOperator{\End}{End}
\DeclareMathOperator{\Endo}{End^{0}}
\DeclareMathOperator{\Gal}{Gal}
\DeclareMathOperator{\Prob}{Prob}
\DeclareMathOperator{\Spec}{Spec}

\newtheorem{theorem}{Theorem}
\newtheorem{lemma}[theorem]{Lemma}

\newtheorem{proposition}[theorem]{Proposition}

\theoremstyle{definition}

\theoremstyle{remark}
\newtheorem{rem}{Remark}

\newtheorem{acks}{Acknowledgments}	
\newtheorem{conventions}{Conventions}	

\begin{document}

\title[Absolutely simple abelian varieties]
{On the existence of absolutely simple\\
abelian varieties of a given dimension\\
over an arbitrary field}
\author{Everett W.\ Howe}
\address{Center for Communications Research, 
         4320 Westerra Court, 
         San Diego, CA 92121-1967, USA.}
\email{however@alumni.caltech.edu}
\urladdr{http://alumni.caltech.edu/\~{}however/}
\author{Hui June Zhu}
\address{Mathematisch Instituut,
         Universiteit Leiden,
         Postbus 9512,
         2300 RA Leiden,
         The Netherlands}
\email{zhu@alum.calberkeley.org}
\date{24 February 2000}


\begin{abstract}
We prove that for every field $k$ and every positive integer~$n$,
there exists an absolutely simple $n$-dimensional abelian variety
over~$k$.  
We also prove an asymptotic result for finite fields:
For every finite field $k$ and positive integer~$n$,
we let $S(k,n)$ denote the fraction of the isogeny classes of
$n$-dimensional abelian varieties over $k$ that consist of absolutely
simple ordinary abelian varieties.
Then for every $n$ we have $S(\BF_{q},n)\ra 1$ as $q\ra\infty$
over the prime powers.
\end{abstract}

\maketitle

\section{Introduction}
\label{S-intro}

An abelian variety over a field $k$ is called {\it simple\/} if 
it has no proper nonzero sub-abelian varieties over $k$; it is called
{\it absolutely simple\/} (or {\it geometrically simple\/})
if it is simple over the algebraic closure of~$k$.
In this paper we will prove the following theorem:
\begin{theorem}
    \label{T-existence}
    Let $k$ be a field and let $n$ be a positive integer.  
    Then there exists an absolutely simple
    $n$-dimensional abelian variety over~$k$.
\end{theorem}

An easy reduction argument, similar in spirit to the one in
Section~4 of~\cite{howe-WH}, shows that an absolutely simple
abelian variety over a field $k$ remains simple over every 
extension field of~$k$, even the non-algebraic ones; thus, it suffices
to prove Theorem~\ref{T-existence} in the special case where $k$ is a
prime field.
Mori~\cite{mori} (see also Zarhin~\cite{zarhin}) provides examples
of absolutely simple abelian varieties of arbitrary dimension over~$\BQ$,
so we need only prove Theorem~\ref{T-existence} for finite prime
fields~$k$.  We will in fact prove that over such fields there exist
absolutely simple {\it ordinary\/} abelian varieties of every dimension,
and in addition we will prove an 
asymptotic result concerning arbitrary finite fields:
\begin{theorem}
    \label{T-asymptotic}
    For every integer $n\ge 0$ and finite field $k$
    let $S(k,n)$ denote the fraction of the isogeny
    classes of $n$-dimensional abelian varieties over $k$ that consist of  
    absolutely simple ordinary abelian varieties.
    Then for every $n$ we have $S(\BF_{q},n)\ra 1$ as $q\ra\infty$
    over the prime powers.
\end{theorem}

In fact, for every $n$ and $\eps$ we provide an explicit value of $M$
such that if $q>M$ then $0 < 1 - S(\BF_{q},n) < \eps$; 
see Theorems~\ref{T-precise-surfaces} and~\ref{T-precise}
in Sections~\ref{S-asymptotic-surfaces} and~\ref{S-asymptotic-higher}.

Suppose $A$ is an abelian variety over a finite field $k$.  One
can ask whether there exists an absolutely simple abelian variety
over $k$ with the same {\it formal isogeny type\/}
(see~\cite{lenstra-oort}) as $A$.
Theorem~\ref{T-existence} shows that the answer to this question
is yes when $A$ is ordinary; we have not considered
the question for other formal isogeny types.
Lenstra and Oort~\cite{lenstra-oort}
considered the analogous question when $k$ is the algebraic closure 
of a finite field, and showed that the answer is yes when $A$ is not 
supersingular.

Theorem~\ref{T-existence} leads to the question of whether there
exist absolutely simple Jacobians of every dimension over a given field~$k$.
Chai and Oort~\cite{chai-oort} show that the answer is yes if 
$k$ is the algebraic closure of a finite field, and Mori~\cite{mori}
and Zarhin~\cite{zarhin} show that the answer is also yes if
$k$ has characteristic zero.  If $k$ has 
positive characteristic $p$ but is not algebraic over~$\BF_{p}$,
then results of Katz and Sarnak (see Sections~10.1 and~10.2 
of~\cite{katz-sarnak}) can be used to show that once again 
the answer is yes ---
see also Mori~\cite{mori} for some partial results for such fields.
In fact, the examples provided by Katz and Sarnak,
Mori, and Zarhin are Jacobians of explicitly-given hyperelliptic curves.
However, the question seems to be open when $k$ is a finite field.
The techniques we use in this paper do not help settle the
general question for finite fields, but our results do at least
show that over every finite field $k$ there are curves of 
genus~$2$ and~$3$ with absolutely simple Jacobians, as the following
argument shows:

As we mentioned above, we show that for every $n$ and for every finite
field $k$ there is an absolutely simple $n$-dimensional ordinary
abelian variety over~$k$, and in particular this is true for
$n=2$ and $n=3$.  But every absolutely simple ordinary abelian
variety of dimension $2$ or $3$ over a finite field is isogenous
to a principally polarized variety (see
Corollary~12.6 and Theorem~1.2 of~\cite{howe-OAV}).
The main result of~\cite{oort-ueno} shows that each such 
principally polarized variety is isomorphic (over the algebraic closure
of~$k$) to a Jacobian of a possibly reducible curve $C$,
but since the Jacobian of $C$ is absolutely simple $C$ must be 
geometrically irreducible.  Finally, a simple descent argument
shows that $C$ has a model defined over~$k$.
Thus, for every finite field $k$ there are curves of genus~$2$ and~$3$
over $k$ with absolutely simple Jacobians.

Our paper is organized as follows.  In Section~\ref{S-Weil}, we 
briefly review the properties of Weil numbers and Weil polynomials
that we will use in the proofs of Theorems~\ref{T-existence}
and~\ref{T-asymptotic}.  In Section~\ref{S-test} we give an
easy-to-verify sufficient condition for an abelian variety over
a finite field to be absolutely simple.  
We use this condition in Section~\ref{S-dimension2} to prove
Theorem~\ref{T-surfaces},
which shows how the characteristic polynomial of Frobenius of a
simple ordinary abelian surface over a finite field can be used 
to quickly determine the splitting behavior of the surface over
the algebraic closure.  Theorem~\ref{T-surfaces} allows us to
give a very short
proof of Theorem~\ref{T-existence} in the case $n=2$; we provide
a proof for the case $n>2$ in Section~\ref{S-higher-dimension}.
In Sections~\ref{S-asymptotic-surfaces} and~\ref{S-asymptotic-higher}
we prove Theorems~\ref{T-precise-surfaces} and~\ref{T-precise},
which are effective versions of Theorem~\ref{T-asymptotic} in
the cases $n=2$ and~$n>2$, respectively.
Finally, in Section~\ref{S-reduction-proof} we prove a lemma about
polynomials with prescribed reduction modulo certain primes that is
essential for our proof of Theorem~\ref{T-precise}.

\begin{conventions}
    Suppose $A$ is an abelian variety over a field $k$ and suppose $\ell$
    is an extension field of~$k$.  We will denote by $A_{\ell}$ the
    $\ell$-scheme~$A\times_{\Spec k} \Spec\ell$.  
    If $B$ is another abelian variety over~$k$, then when we speak of
    a morphism from $A$ to $B$ we always mean a $k$-morphism; 
    thus, we write $\End A$ for what some authors would call~$\End_{k} A$.
\end{conventions}

\begin{acks}
The authors thank 
David Cantor, Robert Coleman,
Daniel Goldstein, Hendrik Lenstra, Bjorn Poonen,
and Joel Rosenberg for helpful conversations and correspondence,
and Mike Zieve for asking the questions that led to this research
and for commenting on an early version of this paper.
The authors are grateful to Hendrik Lenstra for suggesting 
Lemma~\ref{L-subfields}.  
The authors used the computer packages PARI/GP and MAGMA for some of 
the computations they performed in the course of writing this paper.
\end{acks}

\section{Weil numbers and Weil polynomials}
\label{S-Weil}

Suppose $q$ is a power of a prime number~$p$.  A {\it Weil $q$-number},
or simply a {\it Weil number\/} if $q$ is clear from context,
is an algebraic integer $\pi$ such that $|\varphi(\pi)| = q^{1/2}$
for every embedding $\varphi$ of $\BQ(\pi)$ into the complex numbers.
Suppose $k$ is a field with $q$ elements.  
To every abelian variety $A$ over $k$ we associate the characteristic
polynomial $f_{A}\in\BZ[x]$ of its Frobenius endomorphism
(acting on the $\ell$-adic Tate modules of $A$);
the polynomial $f_{A}$ is monic of degree twice the dimension of~$A$.
We call a polynomial $f$ a {\it Weil $q$-polynomial}, or simply a 
{\it Weil polynomial}, if there is an abelian variety $A$ over $k$
with~$f = f_{A}$.
Weil proved that all of the roots of a Weil polynomial are Weil numbers,
and Honda showed that every Weil number is a root of some Weil polynomial.
Furthermore, Tate showed that two abelian varieties over $k$ are
isogenous if and only if their associated Weil polynomials are equal.
If $A$ is a simple abelian variety over $k$ then $f_{A}$ is
a power of an irreducible polynomial, and in fact the Honda-Tate theorem
(see~\cite[Th\'eor\`eme~1]{tate}) says the map that sends $A$ to the
set of roots (in $\BQb$) of $f_{A}$ induces a bijection between
the set of isogeny classes of simple abelian varieties over $k$ and
the set of Galois conjugacy classes of Weil numbers in~$\BQb$.
The Honda-Tate theorem also provides a simple number-theoretic
criterion for determining whether a polynomial, all of
whose roots are Weil numbers, is a Weil polynomial.  In addition,
the theorem shows how the Weil polynomial 
of an abelian variety $A$ over $k$ determines the algebra~$(\End A)\otimes\BQ$.

An abelian variety $A$ over $k$ is {\it ordinary\/} if the rank of
its group of $p$-torsion points over the algebraic closure of $k$ is
equal to the dimension of~$A$; a Weil polynomial is {\it ordinary\/}
if it is the characteristic polynomial of Frobenius of an ordinary 
abelian variety; and a Weil number is {\it ordinary\/} if its
minimal polynomial is an ordinary Weil polynomial.
The Honda-Tate theorem simplifies considerably if one considers only
ordinary varieties and ordinary Weil polynomials --- see
Section~3 of~\cite{howe-OAV}.
For example, a monic polynomial in $\BZ[x]$
is an ordinary Weil $q$-polynomial if and only
if it is of even degree~$2n$, all of its roots are Weil numbers, and
its middle coefficient (that is, the coefficient of $x^{n}$) is
coprime to~$q$.  Furthermore, an ordinary abelian variety $A$ over $k$ is
simple if and only if its Weil polynomial $f$ is irreducible.
If $A$ is simple and ordinary then the algebra $(\End A)\otimes\BQ$
is generated by
the Frobenius endomorphism of~$A$, and so is isomorphic to the
number field defined by~$f$.  Since the characteristic polynomial of 
Frobenius of $A$ has degree equal to twice the dimension of~$A$,
we see that the degree of the number field $K = (\End A)\otimes\BQ$ over 
$\BQ$ is twice the dimension of~$A$.  In fact, the number field $K$
is a {\it CM-field}, which means that 
$K$ is a totally imaginary quadratic extension of a totally real field~$\Kp$.
(A number field $L$ is {\it totally imaginary\/} if it 
cannot be embedded into~$\BR$, and it is 
{\it totally real\/} if every embedding of $L$ into $\BC$
comes from an embedding of $L$ into~$\BR$.)

\section{An easy test for absolute simplicity}
\label{S-test}

In this section we will present an easy-to-verify sufficient
condition for a simple abelian variety over a finite field to be 
absolutely simple.  For ordinary varieties,
the sufficient condition is also necessary.
Throughout this section, $k$ will be a finite field, $\kb$ its 
algebraic closure, $A$ a simple abelian variety over~$k$,
and $\pi$ its Frobenius endomorphism. 
We let $\Endo A$ denote the algebra $(\End A)\otimes\BQ$.
Note that the simplicity of $A$ implies that the subalgebra
$\BQ(\pi)$ of $\Endo A$ is a field.

\begin{proposition}
    \label{P-absolutely}
    Let $D$ be the set of integers $d>1$ such that either 
    \begin{enumerate}
	\item[(a)]
	the minimal polynomial of $\pi$ lies in $\BZ[x^{d}]$ or
	\item[(b)]
	the field $\BQ(\pi^{d})$ is a proper subfield of $\BQ(\pi)$
	and there is a primitive $d$th root of unity $\zeta$ in $\BQ(\pi)$
	such that $\BQ(\pi) = \BQ(\pi^{d},\zeta)$.
    \end{enumerate}
    Then{\rm:}
    \begin{enumerate}
	\item[(1)]
	The set $D$ is empty if and only if
	$\BQ(\pi^{d}) = \BQ(\pi)$ for all $d>0$.
	\item[(2)]
	If $\BQ(\pi^{d}) = \BQ(\pi)$ for all $d>0$ then 
	$A$ is absolutely simple.
	If $A$ is ordinary, then the converse is also true.
    \end{enumerate}
\end{proposition}

To prove this proposition we will need two elementary lemmas.

\begin{lemma}
\label{L-simpleendomorphisms}
    Let $\ell$ be a finite extension of $k$.  
    If $\BQ(\pi^{[\ell : k]}) = \BQ(\pi)$ then $A_{\ell}$ is simple. 
    If $A$ is ordinary, then the converse is also true.
\end{lemma}

\begin{proof}
An abelian variety is simple if and
only if its endomorphism ring contains no zero-divisors.
Thus, if $A$ is simple and $A_{\ell}$ is not, there must
exist an element of $\Endo A_{\ell}$ that does not come from~$\Endo A$.
But it follows from the Honda-Tate theorem~\cite{tate}
that  $\Endo A_{\ell} = \Endo A$ if $\BQ(\pi^{[\ell : k]}) = \BQ(\pi)$.
This proves the first statement of the lemma.

If $A$ is ordinary and $\BQ(\pi^{[\ell : k]})$ is a proper subfield
of $\BQ(\pi)$, then it follows from the Honda-Tate theorem that
$\Endo A_{\ell}$ is a matrix algebra over $\BQ(\pi^{[\ell : k]})$.
In particular, $\Endo A_{\ell}$ contains a zero-divisor, so that 
$A_{\ell}$ is not simple.
\end{proof}

\begin{lemma}
\label{L-subfields}
    Let $\alpha$ be an algebraic number with minimal polynomial $g\in\BQ[x]$,
    and suppose that $d$ is a positive integer such that the field
    $L=\BQ(\alpha^{d})$ is a proper subfield of $K=\BQ(\alpha)$ and such that
    $\BQ(\alpha^{r}) = K$ for all positive~$r < d$.
    Then either $g \in \BQ[x^{d}]$ or there is a primitive $d$th
    root of unity $\zeta$ in $K$ such that~$K=L(\zeta)$.
\end{lemma}

\begin{proof}
Let $\zeta$ be a primitive $d$th root of unity in an algebraic closure
of $K$ and let~$M = L(\zeta)\cap K$.  Note that $M$ contains~$L$.
Since $L(\zeta)$ is a Galois extension of $L$
it is also a Galois extension of~$M$, and it follows that
$L(\zeta)$ and $K$ are linearly disjoint over~$M$, so that
$[K(\zeta) : L(\zeta)] = [K : M]$.
Let $m = [K(\zeta) : L(\zeta)] = [K : M]$.
Since $K(\zeta) = \BQ(\alpha,\zeta)$ is a Kummer extension
of $L(\zeta)=\BQ(\alpha^{d},\zeta)$,
we see that $\alpha^{m}$ lies in~$L(\zeta)$, and hence also in~$M$.

Suppose~$m>1$. Then since $\BQ(\alpha^{m})$ is a subfield of the proper
subfield $M$ of~$K$, the lemma's hypothesis shows we must have~$m=d$.
If we let $h$ be the minimal polynomial of $\alpha^{d}$
over~$\BQ$, then~$g(x) = h(x^{d})$.

Suppose~$m = 1$. Then $K(\zeta) = L(\zeta)$, so 
$K/L$ is a subextension of the abelian 
extension $K(\zeta) / L$, and is therefore Galois.
Let $G$ be its Galois group,
and suppose $\sigma$ is a non-identity element of~$G$.
Let $\xi = \sigma(\alpha)/\alpha$, so that $\xi$ 
lies in the multiplicative group generated by~$\zeta$.
Suppose $r$ is a positive integer less than~$d$.  
Then the hypothesis of the lemma shows that
$K= \BQ(\alpha^{r})$, so we must have 
$\alpha^{r}\neq\sigma(\alpha^{r}) = \xi^{r}\alpha^{r}$.
Thus $\xi$ must in fact be a primitive $d$th root of unity, which 
shows that~$\zeta\in K$.  It follows that 
$K= K(\zeta)$, and this last field is $L(\zeta)$ 
because~$m = 1$.
\end{proof}

\begin{proof}[Proof of Proposition~{\rm\ref{P-absolutely}}]
If $d$ is an integer in $D$ then clearly $\BQ(\pi^{d})$ is a proper
subfield of $\BQ(\pi)$.  On the other hand, if there exists some $d>0$
such that $\BQ(\pi^{d})\neq\BQ(\pi)$ then there exists a smallest
such~$d$, and by Lemma~\ref{L-subfields} this $d$ lies in $D$.
This proves the first statement of the proposition.

It is clear that $A$ is absolutely simple if and only if $A_{\ell}$ is
simple for every finite extension $\ell$ of~$k$.  The second
statement of the proposition follows from this fact 
and from Lemma~\ref{L-simpleendomorphisms}.
\end{proof}

\begin{rem}
A theorem of Silverberg \cite{silverberg} shows that 
if $A$ is an abelian variety over an arbitrary field~$k$,
then to check that 
$\Endo A = \Endo A_{\kb}$ it suffices to check
that $\Endo A = \Endo A_{\ell}$
for a certain finite extension $\ell$ of $k$;
in particular, if one chooses an integer $m>2$ not divisible
by the characteristic of~$k$, Silverberg shows that
one may take $\ell$ to be the smallest field
over which every $m$-torsion point of $A$ is defined.
The degrees of such $\ell$ over $k$ may be quite large, even when
$k$ is a finite field.
Lemmas~\ref{L-subfields} and the proof of Lemma~\ref{L-simpleendomorphisms}
show that Silverberg's general result can be improved in the
special case where $k$ is finite.
\end{rem}

\section{Absolutely simple abelian surfaces}
\label{S-dimension2}

In this section we will prove a theorem that shows that,
given the characteristic polynomial of Frobenius of a simple
ordinary abelian surface over a finite field, it is quite easy to
determine whether the surface is absolutely simple.
At the end of the section we will use this theorem to prove
the special case $n=2$ of Theorem~\ref{T-existence}.

Suppose $k$ is a finite field with $q$ elements and $A$ is an
abelian surface over~$k$.  Let $f$ be the characteristic
polynomial of Frobenius for~$A$.  Then Weil's ``Riemann Hypothesis''
shows that $f$ is of the form $x^{4}+ax^{3}+bx^{2}+aqx+q^{2}$
for some integers $a$ and~$b$.  
If neither $a$ nor $b$ is coprime to $q$ then one can
use the Honda-Tate theorem to show that $A$ becomes isogenous to
the square of a supersingular elliptic curve over a finite extension
of~$k$.  If $a$ is coprime to $q$ but $b$ is not, then one can
again use Honda-Tate to show that $A$ is absolutely simple if and
only if it is simple, and that $A$ is simple if and only if
$f$ is irreducible.  
The most interesting situation arises when $b$ is coprime to~$q$,
which is the case exactly when $A$ is an ordinary abelian variety.
In this case, $A$ is simple if and only if $f$ is irreducible.

\begin{theorem}
    \label{T-surfaces}
    Suppose $f = x^{4}+ax^{3}+bx^{2}+aqx+q^{2}$ is the Weil polynomial
    of a simple ordinary abelian surface $A$ over~$k$.  Then exactly one
    of the following conditions holds{\rm:}
    \begin{enumerate}
	\item[(a)] The variety $A$ is absolutely simple.
	\item[(b)] We have $a=0$.
	\item[(c)] We have $a^{2} = q + b$.
	\item[(d)] We have $a^{2} = 2b$.
	\item[(e)] We have $a^{2} = 3b - 3q$.
    \end{enumerate}
    In cases \text{\rm(b)}, \text{\rm(c)}, \text{\rm(d)}, and \text{\rm(e)},
    the smallest extension of $k$ over which $A$ splits is
    quadratic, cubic, quartic, and sextic, respectively.
\end{theorem}

\begin{proof}
Let $\pi$ be the Frobenius endomorphism of $A$ and let $K$ be the 
field~$\BQ(\pi)$.  Because $A$ is ordinary, the field $K$ is a
CM-field of degree $4$ over~$\BQ$.
The ordinariness of $A$ also implies that $\BQ(\pi^{d})$
is a CM-field for every positive integer~$d$, and that $A$ splits
over the degree-$d$ extension of $k$ if and only if $\BQ(\pi^{d})$
is a proper subfield of~$K$.

Suppose $A$ is not absolutely simple.  Then there is a positive integer $d$
such that $\BQ(\pi^{d})$ is a proper subfield of $K$; let us take $d$ to be
the smallest such integer, and let $L$ be the imaginary quadratic 
field~$\BQ(\pi^{d})$.  
By Lemma~\ref{L-subfields},
either $d = 2$ and $f\in\BZ[x^{2}]$, or $d=4$ and $f\in\BZ[x^{4}]$,
or there is a primitive $d$th root of unity $\zeta$ in $K$ such 
that~$K = L(\zeta)$.  Let us first show that if the third possibility
is the case and if $d>4$ then $d$ must equal~$6$.

Suppose, to obtain a contradiction, that we are in the third case
and that $d$ is greater than $4$ but not equal to~$6$.
Then the degree of $\BQ(\zeta)$ over $\BQ$ is greater than~$2$,
so $K$ must be~$\BQ(\zeta)$.
Let $\sigma\in\Gal(K/\BQ)$ be the nontrivial automorphism of $K$ that
fixes~$L$.  The proof of Lemma~\ref{L-subfields} shows that 
we may choose our primitive root of unity $\zeta$ so
that~$\pi^{\sigma} = \zeta\pi$.  Applying $\sigma$ to this equality,
we find $\pi = \zeta^{\sigma}\pi^{\sigma} = \zeta^{\sigma}\zeta\pi$,
so that~$\zeta^{\sigma}\zeta = 1$.  The only element of the Galois group
of the cyclotomic field with this property is complex conjugation.
But then the fixed field $L$ of $\sigma$ must be totally real, and
we have reached a contradiction.

So we must find, for $d = 2, 3, 4,$ and $6$, the conditions on the
coefficients $a$ and $b$ in the minimal polynomial of $\pi$ that
are equivalent to $\pi^{d}$ lying in a quadratic subfield.  
Note that the characteristic polynomial of $\pi^{d}$ is 
of the form $x^{4}+\alpha x^{3}+\beta x^{2}+\alpha q^{d}x+q^{2d}$,
and that such a quartic polynomial is the square of a quadratic
polynomial if and only if $\alpha^{2} - 4\beta + 8q^{d} = 0$.
It is not difficult to explicitly calculate the characteristic 
polynomial of $\pi^{d}$ for each $d$ we are considering, and we find
that
\begin{equation*}
    \alpha^{2} - 4\beta + 8q^{d} =
    \begin{cases}
	a^{2} (a^{2}-4b + 8q)           & \text{if $d=2$;}\\
        (a^{2}-b-q)^{2} (a^{2}-4b + 8q)     & \text{if $d=3$;}\\
        a^{2} (a^{2} - 2b)^{2} (a^{2}-4b + 8q) & \text{if $d=4$;}\\
	a^{2} (a^{2}-b-q)^{2} (a^{2} - 3b +3q)^{2} (a^{2}-4b + 8q) &
	               \text{if $d=6$.}
    \end{cases}
\end{equation*}
We have assumed that $A$ is simple over~$k$, so the characteristic 
polynomial for $\pi$ is irreducible; this means in particular
that the quantity $a^{2} - 4b + 8q$ is nonzero.  Thus, if $A$ is
not absolutely simple then one of the cases (b), (c), (d), or (e) 
must hold.  Note that if two of these cases
were to hold simultaneously, then $b$ would equal a multiple of~$q$,
contradicting our assumption that $b$ is coprime to~$q$. Thus
exactly one of the cases (a) through (e) must hold.

Finally, the formulas for $\alpha^{2} - 4\beta + 8q^{d}$ given above
make it easy to verify the theorem's statement about 
the degree of the minimal splitting field of~$A$.
\end{proof}

Using Theorem~\ref{T-surfaces}, it is easy to show that there exist
absolutely simple ordinary abelian surfaces over every finite field.
If $q$ is an arbitrary prime power, then
Theorem~1.1 of~\cite{ruck} shows that the polynomial
$x^{4} + x^{3} + x^{2} + qx + q^{2}$ is an ordinary Weil
polynomial.  It is easy to check that this polynomial is
irreducible, so it corresponds to an isogeny class of simple
abelian varieties over the field~$\BF_{q}$.  
Then Theorem~\ref{T-surfaces} shows
that the varieties in the isogeny class are absolutely simple.

\section{The existence of absolutely simple abelian varieties of higher dimension}
\label{S-higher-dimension}

In this section we will prove Theorem~\ref{T-existence}
in the case where~$n>2$.  As we noted in
the Introduction, it suffices to prove the theorem for finite
prime fields~$k$, but we will assume only that $k$ is finite.  
In fact, for such fields we will prove a result that
is slightly stronger than Theorem~\ref{T-existence}.

\begin{theorem}
    \label{T-higher}
    Let $k$ be a finite field and let $n>2$ be an integer.
    Then there is an absolutely simple $n$-dimensional ordinary abelian 
    variety over~$k$.
\end{theorem}

The proof of Theorem~\ref{T-higher} depends on three lemmas,
whose proofs we will postpone until after the proof of the theorem.
The first lemma gives sufficient conditions for an ordinary Weil
number to correspond to an isogeny class of absolutely simple
varieties.

\begin{lemma}
    \label{L-good-Weil-numbers}
    Let $q$ be a prime power and let $n>2$ be an integer. 
    Suppose $\pi$ is an ordinary Weil $q$-number, let $K=\BQ(\pi)$, 
    let $\Kp$ be the maximal real subfield of~$K$, and let~$n = [\Kp:\BQ]$. 
    Suppose that
    \begin{enumerate}
	\item[(1)] the minimal polynomial of $\pi$ is not of 
	the form $x^{2n} + ax^n + q^n$, 
	\item[(2)] the field $\Kp$ has no proper subfields
	other than~$\BQ$, and
	\item[(3)] the field $\Kp$ is not the maximal real
	subfield of a cyclotomic field.
    \end{enumerate}
    Then the isogeny class corresponding to $\pi$ consists of 
    absolutely simple varieties.
\end{lemma}

The second lemma shows that any polynomial satisfying a certain
set of local conditions also satisfies the
hypotheses of Lemma~\ref{L-good-Weil-numbers}.
We will use this lemma again in Section~\ref{S-asymptotic-higher}.

\begin{lemma}
    \label{L-good-polynomials}
    Let $q$ be a prime power and let $n>2$ be an integer. 
    Let $g\in\BZ[x]$ be a monic polynomial of degree~$n$, and let $f$ be the
    polynomial given by $f(x) = x^{n} g(x+q/x)$.  Suppose that
    the following five conditions hold\/{\rm{:}}
    \begin{enumerate}
	\item[(1)]  the polynomial $f$ is not of 
	the form $x^{2n} + ax^n + q^n$, 
	\item[(2)]  all of the complex roots of $g$ are real numbers of 
	absolute value less than~$2\sqrt{q}$,
	\item[(3)]  the constant term of $g$ is coprime to~$q$,
	\item[(4)]  there exists a prime $p_{1}$ such that the
	reduction of $g$ modulo $p_{1}$ is irreducible, and
	\item[(5)]  there exists a prime $p_{2}$ such that the 
	reduction of $g$ modulo $p_{2}$ is a linear times 
	an irreducible.
    \end{enumerate}
    Then $f$ is an irreducible ordinary Weil polynomial of degree~$2n$,
    and its roots $\pi$ satisfy the hypotheses of 
    Lemma~{\rm\ref{L-good-Weil-numbers}}.
\end{lemma}

The third lemma gives us a way of producing polynomials that meet
the hypotheses of Lemma~\ref{L-good-polynomials}.

\begin{lemma}
    \label{L-good-polys-exist}
    Let $q$ be a prime power and let $n>2$ be an integer. Then there is
    a monic polynomial $g$ in $\BZ[x]$ that satisfies the following
    five conditions\/{\rm{:}}
    \begin{enumerate}
	\item[(1)]  the polynomial $g$ can be written
        $$g = x^n + cx^{n-2} + \text{lower-order terms},$$
	where either $c$ is equal to $-2n$ or $c$ is not divisible by~$n$,
	\item[(2)]  all of the complex roots of $g$ are real numbers of 
	absolute value less than~$2\sqrt{2}$, 
	\item[(3)]  the constant term of $g$ is coprime to~$q$,
	\item[(4)]  the reduction of $g$ modulo $2$ is irreducible, and
	\item[(5)]  the reduction of $g$ modulo $3$ is a linear times 
	an irreducible.
    \end{enumerate}
\end{lemma}
    
\begin{proof}[Proof of Theorem~{\rm\ref{T-higher}}]
Let $g$ be the polynomial whose existence is guaranteed by 
Lemma~\ref{L-good-polys-exist}.  Then $g$ satisfies the last four of the 
five hypotheses of Lemma~\ref{L-good-polynomials};
we will show that it satisfies the first hypothesis as well.

First we will consider the case in which~$q>2$.  Since
      $$g = x^n + cx^{n-2} + \text{lower-order terms},$$
we find that the polynomial $f$ defined in Lemma~\ref{L-good-polynomials}
may be written in the form
      $$f = x^{2n} + (qn + c)x^{2n-2} + \text{lower-order terms}.$$
Now, $c$ is either $-2n$ or is not a multiple of~$n$, so the coefficient
of $x^{2n-2}$ in $f$ is nonzero.  In particular, 
$f$ is not of the form $x^{2n} + ax^n + q^n$.

For the case in which $q = 2$ we use the easily-proven fact that the 
reduction of $f$
modulo $2$ is equal to $x^n$ times the reduction of $g$ modulo~$2$.
Since $g$ modulo $2$ is irreducible, and since $x^n + 1$ is not irreducible
over~$\BF_2$, the polynomial $f$ must have an odd coefficient somewhere
between $x^{2n}$ and~$x^n$.  Again we see that $f$ is not of 
the form $x^{2n} + ax^n + q^n$.

Thus $g$ satisfies all the hypotheses of Lemma~\ref{L-good-polynomials},
so by Lemma~\ref{L-good-Weil-numbers} the roots of $f$ are
Weil numbers that correspond to an isogeny class of absolutely simple 
ordinary varieties over~$\BF_{q}$.
\end{proof}

\begin{proof}[Proof of Lemma~{\rm\ref{L-good-Weil-numbers}}]
Suppose, to obtain a contradiction,
that $\pi$ corresponds to an isogeny class that is not absolutely
simple. Then by Proposition~\ref{P-absolutely}
there is a positive integer $d$ such that
$\BQ(\pi^{d})$ is a proper subfield of~$K$.  Let $d$ be the
smallest positive integer with this property.  Since $\pi$ is
ordinary, the field $L = \BQ(\pi^{d})$ is a CM-field, and its maximal real
subfield $L^{+}$ is a proper subfield of~$\Kp$.
Hypothesis (2) shows that $L^{+}$ must be~$\BQ$, so $L$ is an imaginary
quadratic field.

Lemma~\ref{L-subfields} shows that either the minimal polynomial $f$ of $\pi$
lies in $\BZ[x^{d}]$ or $K = L(\zeta)$ for some primitive $d$th root of
unity. 
The first possibility cannot occur, because it would imply that~$d=n$,
contradicting hypothesis (1).
Therefore the second possibility must be
the case.  We find that the maximal real subfield of $\BQ(\zeta)$ is 
a subfield of~$\Kp$, and since $\Kp$ is not itself the maximal real subfield of
a cyclotomic field  (by assumption), we find that the maximal real subfield
of $\BQ(\zeta)$ must be~$\BQ$, so that $\BQ(\zeta)$ is either a quadratic
field or $\BQ$ itself.  But $K$ is the compositum of $L$ and
$\BQ(\zeta)$, so the degree of $K$ over $\BQ$ is at most~$4$.
This contradicts our assumption that the degree of $K$ over $\BQ$ is~$2n$, 
where~$n>2.$
\end{proof}

\begin{proof}[Proof of Lemma~{\rm\ref{L-good-polynomials}}]
Since $g$ modulo $p_{1}$ is irreducible, $g$ itself is irreducible in $\BZ[x]$,
and since all of its complex roots are real, $g$ defines a totally real number
field~$\Kp$.  Let $\alpha$ be a root of $g$ in~$\Kp$.  The discriminant
of the polynomial $h = x^{2} - \alpha x + q$ is totally negative because the
roots of $g$ all have magnitude less than~$2\sqrt{q}$, so $h$ defines
a totally imaginary quadratic extension $K$ of~$\Kp$.  If $\pi$ is
a root of $h$ in~$K$, then $K=\BQ(\pi)$ contains $\Kp$ because 
$\alpha = \pi + q/\pi$.  Thus $\pi$ is an algebraic number of degree~$2n$.
Furthermore, if $\varphi$ is an embedding of $K$ into~$\BC$, then
$\varphi(\pi)$ is a root of $x^{2} - \varphi(\alpha)x + q$, and 
the quadratic formula shows that $|\varphi(\pi)| = \sqrt{q}$.
Thus $\pi$ is in fact a Weil number of degree~$2n$.
Since $\pi$ is a root of~$f$, the polynomial $f$ must be the minimal polynomial
of~$\pi$.  This shows that $f$ is an irreducible polynomial whose roots
are Weil numbers, and to show that $f$ is an ordinary Weil polynomial
we need merely check that its middle coefficient is coprime to~$q$.
But this follows from hypothesis (3), because the middle coefficient of $f$
differs from
the constant term of $g$ by a multiple of~$q$.

Now we must check that a root $\pi$ of $f$ satisfies the hypotheses of
Lemma~\ref{L-good-Weil-numbers}.  The first of these hypotheses is identical to
the first hypothesis of the lemma we are proving, and is therefore satisfied.

We will show that $\Kp$ is not the maximal real subfield of 
a cyclotomic field.  It will suffice to show that $\Kp$
is not Galois over~$\BQ$.  The defining polynomial $g$
of $\Kp$ reduces modulo $p_{2}$ as a linear times an irreducible, so the
prime $p_{2}$ splits in $\Kp$ into two primes with different residue class 
degrees, so $\Kp/\BQ$ cannot be Galois.

Finally, we prove that $\Kp$ has no proper subfields other than~$\BQ$.
For suppose $\Kp$ had a proper subfield $L$ other than~$\BQ$.  Let $\Pid$ be
the prime  of $\Kp$ over $p_{2}$ whose residue class degree is~$n-1$.
Let $\pid$ be the prime of $L$ lying under~$\Pid$.  Let $f_{1}$ be the residue 
class degree of $\Pid$ over $\pid$ and let $f_{2}$ be the residue class degree
of $\pid$ over~$p_{2}$.  Then we have the three statements:
\begin{enumerate}
    \item[(a)] $f_{1} \le [\Kp : L]$,
    \item[(b)] $f_{2} \le [L : \BQ]$, and
    \item[(c)] $f_{1} f_{2} = n - 1 = -1 + [\Kp : L] [L : \BQ].$
\end{enumerate}
Statement (c) shows that strict inequality must hold in one
of statements (a) and~(b); but then, since both of the field
extensions $\Kp/L$ and $L/\BQ$ are assumed non-trivial, we find that 
$f_{1} f_{2}$ must be {\it less\/} than $-1 + [\Kp : L] [L : \BQ]$. 
This contradiction shows that $\Kp$  has no proper subfields other than~$\BQ$.
\end{proof}

Our proof of Lemma~\ref{L-good-polys-exist} depends on a result of
Robinson concerning certain modified Chebyshev polynomials.
Before starting on the proof of the lemma
we will define these polynomials and present Robinson's result.

For every positive integer $i$ let $t_i$ be the $i$th Chebyshev polynomial,
so that $t_i(x) = \cos(i \cdot\arccos(x))$.  For every positive integer $i$
let $T_i$ be the polynomial given by 
$T_i(x) = 2\cdot 2^{i/2} \cdot t_i(x/2^{3/2})$.
It is not hard to show that $T_i$
is a monic polynomial in~$\BZ[x]$
and that $T_{i}\equiv x^{i} \bmod 2$.  Let~$T_0 = 1$.

\begin{lemma}
    \label{L-Robinson}
    Suppose $a_{1},\ldots,a_{n}$ are real numbers such that
    $$\left(\sum_{i=1}^{n-1} \left| \frac{a_{i}}{2^{i/2}} \right|\right)
       + \frac{1}{2}\left|\frac{a_{n}}{2^{n/2}}\right| < 1.$$
    Then every complex root of the polynomial
    $$T_{n} + a_{1}T_{n-1} + \cdots + a_{n-1}T_{1} + a_{n}T_{0}$$
    is real and lies in the open interval $(-2\sqrt{2},2\sqrt{2})$.
\end{lemma}

\begin{proof}
This follows from the techniques of Robinson~\cite{robinson}.
\end{proof}

\begin{proof}[Proof of Lemma~{\rm\ref{L-good-polys-exist}}]
If $n\le 9$ we can simply choose the
appropriate value of $g$ from Table~\ref{table-smallgood},
so let us assume that~$n>9$.

\begin{table}
\renewcommand\arraystretch{1.25}
\begin{center}
\begin{tabular}{|c|c|}
\hline
$n$ & $g$                                                     \\ \hline\hline
$3$ & $x^3 - 5x + 1$                                          \\ \hline
$4$ & $x^4 - 6x^2 - x + 1$                                    \\ \hline
$5$ & $x^5 - 10x^3 + x^2 + 20x + 1$                           \\ \hline
$6$ & $x^6 - 12x^4 + 34x^2 + x - 1$                           \\ \hline
$7$ & $x^7 - 14x^5 + 56x^3 - 2x^2 - 57x + 1$                  \\ \hline
$8$ & $x^8 - 16x^6 + 81x^4 + x^3 - 129x^2 + 1$                \\ \hline
$9$ & $x^9 - 18x^7 + 108x^5 + x^4 - 240x^3 - 9x^2 + 147x + 1$ \\ \hline
\end{tabular}
\end{center}
\vspace{1ex}
\caption{Polynomials satisfying the conditions of Lemma~\ref{L-good-polys-exist}
for small values of~$n$.}
\label{table-smallgood}
\end{table}

Lemma~\ref{L-mod-2-and-3} (below) shows that there exist monic 
degree-$n$ polynomials $g_2$
in $\BF_2[x]$ and $g_3$ in $\BF_3[x]$ such that $g_2$ is irreducible, 
such that $g_3$ is a linear times an irreducible and has nonzero constant term,
and such that the coefficients of $x^{n-1}, \ldots, x^{n-6}$ in $g_2$ and $g_3$
are equal to the reductions (modulo $2$ and $3$) of the corresponding
coefficients of the modified Chebyshev polynomial~$T_n$.
Once we have fixed $g_{2}$ and~$g_{3}$,
we can choose values of $a_7, a_8, \ldots, a_n$ in the set $\{-2,-1,0,1,2,3\}$
such that the polynomial
$$  g = T_n + a_7 T_{n-7} + a_8 T_{n-8} + ... + a_{n-1} T_1 + a_n T_0  $$
reduces to $g_2$ modulo $2$ and to $g_3$ modulo~$3$.

Note that the constant term of $g$ is coprime to $6$ because $g_2$ and $g_3$
have nonzero constant terms.  Thus, if $q$ is a power of $2$ or $3$ then
the constant term of $g$ is coprime to~$q$.  If $q$ is not a power of $2$ or~$3$,
then the constant term of $g$ may have a factor in common with~$q$.  
If this is the case, replace $a_n$ with either $a_n - 6$ or $a_n + 6$, 
whichever one lies in the interval $[-6, 6]$; this changes the constant
term of $g$ by~$6$, so that the constant term is now coprime to $q$ 
but so that $g$ still reduces to $g_2$ modulo $2$ and to $g_3$ modulo~$3$.

One can calculate that 
$T_{n} = x^n - 2n x^{n-2} + \text{lower-order terms},$
and since $g$ differs from $T_{n}$ by a polynomial of degree at
most $n-7$, we see that $g$ may also be written 
$g = x^n - 2n x^{n-2} + \text{lower-order terms}.$
In particular, $g$ satisfies the first condition of 
Lemma~\ref{L-good-polys-exist}.

Thus $g$ satisfies four of the five conditions listed in
the statement of Lemma~\ref{L-good-polys-exist}.
We are left to show that all of its roots are real, and that they
have absolute value less than~$2\sqrt{2}$.
But this follows from Lemma~\ref{L-Robinson};
to apply the lemma we must verify that
the quantity 
$$\left| \frac{a_{7}}{2^{7/2}} \right|
 + \left| \frac{a_{8}}{2^{8/2}} \right|
 + \cdots 
 + \left| \frac{a_{n-1}}{2^{(n-1)/2}} \right|
 + \frac{1}{2} \left| \frac{a_{n}}{2^{n/2}} \right|$$
is less than~$1$, and this follows from the fact that $|a_{i}|$ is at 
most $3$ for~$i<n$, and that $|a_{n}|$ is at most~$6$.

Thus, the $g$ we have written down satisfies all the conditions of 
the lemma.
\end{proof}

\begin{lemma}
    \label{L-mod-2-and-3}
    Suppose $n\ge 10$.  Then there exist monic degree-$n$
    polynomials $g_2$ in $\BF_2[x]$ and $g_3$ in $\BF_3[x]$ such that 
    $g_2$ is irreducible, such that 
    $g_3$ is a linear times an irreducible and has
    nonzero constant term, and such that
    the coefficients of $x^{n-1}$ through $x^{n-6}$ of $g_2$ and $g_3$
    are equal to the reductions modulo $2$ and $3$ of the corresponding
    coefficients of the modified Chebyshev polynomial $T_n$
    defined above.
\end{lemma}

\begin{proof}
For $n \le 18$ we choose $g_{2}$ and $g_{3}$ from
Table~\ref{table-2-and-3}.  For $n > 18$ we argue as follows:

\begin{table}
\renewcommand\arraystretch{1.25}
\begin{center}
\begin{tabular}{|c|c|c|}
\hline
$n$ &      $g_2$          &                 $g_3$                 \\ \hline\hline
$10$&$x^{10}+x^3+1$       &$x^{10}+x^8   -x^6   -x^4   +x^2+x  +1$\\ \hline
$11$&$x^{11}+x^2+1$       &$x^{11}-x^9   -x^7   -x^5   +1$        \\ \hline
$12$&$x^{12}+x^3+1$       &$x^{12}+x^6   -x^2   +x     +1$        \\ \hline
$13$&$x^{13}+x^5+x^2+x+1$ &$x^{13}+x^{11}-x^9   +x^2   +1$        \\ \hline
$14$&$x^{14}+x^5+1$       &$x^{14}-x^{12}-x^{10}+x^2   +x  +1$    \\ \hline
$15$&$x^{15}+x  +1$       &$x^{15}-x^9   +x     +1$               \\ \hline
$16$&$x^{16}+x^6+x^2+x+1$ &$x^{16}+x^{14}-x^{12}+x^{10}+x^2+x  -1$\\ \hline
$17$&$x^{17}+x^3+1$       &$x^{17}-x^{15}-x^{13}+x^{11}+x^3+x^2+1$\\ \hline
$18$&$x^{18}+x^3+1$       &$x^{18}+x^2   +x     -1$               \\ \hline
\end{tabular}
\end{center}
\vspace{1ex}
\caption{Polynomials satisfying the conditions of Lemma~\ref{L-mod-2-and-3}
for small values of~$n$.}
\label{table-2-and-3}
\end{table}

Corollary 3.2 (p. 94) of~\cite{hsu} shows
that there exists a monic irreducible polynomial in $\BF_2[x]$ of degree $n$
with zeroes for the first six coefficients after the leading~$x^n$.
We take this polynomial for our~$g_2$.  The same corollary shows that
there is a monic irreducible polynomial $h$ in $\BF_3[x]$ such that the
first six coefficients of $(x-1)h$ are equal to those of the 
reduction of $T_n$ modulo $3$; we take $g_3$ to be~$(x-1)h$.
\end{proof}

\section{Asymptotic results for abelian surfaces}
\label{S-asymptotic-surfaces}

In this section we will prove Theorem~\ref{T-asymptotic} 
in the case~$n = 2$.  In fact, we will prove a more precise statement.

\begin{theorem}
    \label{T-precise-surfaces}
    Let $\eps$ be a positive real number.
    If $q$ is a prime power with $q > (659/\eps)^{2}$ then
    $S(\BF_{q},2) > 1-\eps$.
\end{theorem}

\begin{proof}
Let $r$ be the arithmetic function defined by $r(x) = \varphi(x)/x$, 
where $\varphi$ is Euler's $\varphi$-function,
let $I$ be the number of isogeny classes
of abelian surfaces over~$\BF_{q}$, let $O_{\simp}$ be the number of isogeny
classes of simple ordinary abelian surfaces, and let $O_{\asimp}$ 
be the number of isogeny classes of absolutely simple ordinary abelian 
surfaces.  Theorem~1.2
of~\cite{dipippo-howe} shows that
$$I < \frac{32}{3} r(q) q^{3/2} + 3473 q + 8359 q^{1/2};$$
this upper bound is obtained by combining the estimates that
Theorem~1.2 gives for the number of ordinary and non-ordinary
isogeny classes of abelian surfaces.  

The same theorem shows that the number of isogeny classes of ordinary 
abelian surfaces over $\BF_{q}$ is at least
$$\frac{32}{3} r(q) q^{3/2} - 8359 q^{1/2}.$$

The isogeny classes of ordinary elliptic curves over $\BF_{q}$ correspond
to the integers $t$ such that $|t| < 2q^{1/2}$ and $(t,q) = 1$, so there
are at most $4q^{1/2}$ such isogeny classes.  A non-simple isogeny class
of ordinary abelian surfaces is determined by its two factors,
so there are at most $4q^{1/2}(4q^{1/2} + 1)/2 = 8q + 2q^{1/2}$ such 
reducible isogeny classes.  Thus we have 
$$O_{\simp} > \frac{32}{3} r(q) q^{3/2} - 8q - 8361 q^{1/2}.$$

Now we must estimate the number of simple ordinary isogeny classes that are not
absolutely simple.  For this we use Theorem~\ref{T-surfaces}.
First note that if $x^{4} + ax^{3} + bx^{2} + aqx + q^{2}$ is 
the Weil polynomial for an ordinary abelian surface over $\BF_{q}$ then
$|a| < 4q^{1/2}$, and if $a = 0$ then $0 < |b| < 2q.$
Thus, the number of Weil polynomials of ordinary abelian surfaces
that satisfy case (b) of Theorem~\ref{T-surfaces} is at most~$4q$.
Also, for every nonzero integer $d$ in the interval $(-4q^{1/2}, 4q^{1/2})$ there
is at most one Weil polynomial with $a = d$ that satisfies 
case (c) of the theorem; for every nonzero integer $d$ in the interval
$(-2q^{1/2}, 2q^{1/2})$ there
is at most one Weil polynomial with $a = 2d$
that satisfies case (d) of the theorem; and for 
every nonzero  integer $d$ in the interval
$(-(4/3)q^{1/2}, (4/3)q^{1/2})$ there
is at most one Weil polynomial with $a = 3d$ 
that satisfies case (e) of the theorem.
We find that there are at most $15q^{1/2}$ simple Weil polynomials
$x^{4} + ax^{3} + bx^{2} + aqx + q^{2}$ with $a\neq 0$ that
are not absolutely simple.  

Combining these estimates with the lower bound for $O_{\simp}$ given above,
we find that 
$$O_{\asimp} >  \frac{32}{3} r(q) q^{3/2} - 12 q - 8376 q^{1/2}.$$

Now suppose $\eps$ is given.  If $\eps\ge 1$ then the conclusion of the
theorem is clearly true for all~$q$, so we may assume that $\eps < 1$
and that~$q>659^{2}$.  With this lower bound for~$q$, our bounds for $I$
and $O_{\asimp}$ show that 
$$ I < \frac{32}{3} r(q) q^{3/2} + 3486 q $$
and
$$ O_{\asimp} > \frac{32}{3} r(q) q^{3/2} - 25 q.$$
Thus
$$\frac{O_{\asimp}}{I} > \left({1 - \frac{75}{32 r(q) q^{1/2}}}\right)
                \bigg/
		\left({1 + \frac{10458}{32r(q) q^{1/2}}}\right).$$
The denominator is less than~$2$, so we have
\begin{align*}
\frac{O_{\asimp}}{I} &>\left(1 - \frac{75}{32 r(q) q^{1/2}}\right)
                        \left(1 - \frac{10458}{32r(q) q^{1/2}}\right)\\
	             &>1 - \frac{10533}{32 r(q) q^{1/2}}\\
	             &>1 - \frac{10533}{16q^{1/2}}
	              >1 - \frac{659}{q^{1/2}}
		      >1 - \eps,
\end{align*}
as was to be shown.
\end{proof}

\section{Asymptotic results for abelian varieties of higher dimension}
\label{S-asymptotic-higher}

In this section we will prove Theorem~\ref{T-asymptotic} 
in the case $n > 2$ by 
proving a more precise result, whose statement requires
us to introduce some notation. 
First we define constants $c_{1}, c_{2},$ and $c_{3}$ by setting
\begin{align*}
c_1 &= \sqrt{3}/6\approx 0.288675, \\
c_2 &= \exp(3/2) \cdot 2 (1+\sqrt{2}) \sqrt{3} (1+\sqrt{3}/162)^3/3   
       \approx 12.898608,\\
\intertext{and}
c_3 &= c_2 / (1+\sqrt{2}) \approx 5.342778.
\end{align*}
Next, for every positive integer $n$ we let 
$$v_{n}= \frac{2^{n}}{n!} \prod_{j=1}^{n} \left(\frac{2j}{2j-1}\right)^{n+1-j}$$
and we let 
$$G_{n} = \frac{1}{v_{n}} 6^{n^{2}} c_{1}^{n} c_{3} \frac{n(n+1)}{(n-1)!}.$$
Finally, if $n>1$ is an integer and if $\eps$ is a positive real,
we let $k_{n,\eps}$ denote the smallest positive integer $k$ such that
$$\left(1-\frac{1}{2n}\right)^{k}  \le \frac{\eps}{8},$$
we let $m_{n,\eps}$ be the product of
the first $k_{n,\eps}$ prime numbers, and we let 
$$M_{n,\eps} = \left(\frac{8 G_{n} m_{n,\eps}}{\eps}\right)^{2}.$$
Recall that $S(\BF_{q},n)$ denotes the fraction of isogeny classes
of $n$-dimensional abelian varieties over $\BF_{q}$ that are ordinary
and absolutely simple.

\begin{theorem}
    \label{T-precise}
    Let $n>2$ be an integer and let $\eps$ be a positive
    real number. If $q > M_{n,\eps}$ then $S(\BF_{q},n) > 1-\eps.$
\end{theorem}

For every prime power $q$ and non-negative integer $n$ we let $\CI(q,n)$
denote the set of isogeny classes of $n$-dimensional abelian varieties
over $\BF_{q}$ and we let $\CO(q,n)$ and $\CN(q,n)$ denote
the sets of ordinary and non-ordinary isogeny classes in~$\CI(q,n)$,
respectively.
Also, we let $\CO_{\simp}(q,n)$ and $\CO_{\asimp}(q,n)$ denote the
sets of simple and absolutely simple isogeny classes in~$\CO(q,n)$,
respectively.
As in Section~\ref{S-asymptotic-surfaces}
we let $r$ be the arithmetic function
defined by $r(x) = \varphi(x)/x$, where $\varphi$ is Euler's $\varphi$-function.
Our proof of Theorem~\ref{T-precise} breaks into two parts.  First we
will give an upper bound for~$\#\CI(q,n)$.

\begin{proposition}
    \label{P-upper-I}
    Let $n>2$ be an integer and let $\eps$ be a positive
    real number with $\eps \le 1$. If $q > M_{n,\eps}$ then 
    $\#\CI(q,n) < (1 + \eps/8) v_{n} r(q) q^{n(n+1)/4}.$
\end{proposition}

Then we will give  a lower bound for $\#\CO_{\asimp}(q,n)$.

\begin{proposition}
    \label{P-lower-O}
    Let $n>2$ be an integer and let $\eps$ be a positive
    real number with $\eps\le 1$. If $q > M_{n,\eps}$ then 
    $\#\CO_{\asimp}(q,n) \ge (1 - 7\eps/8) v_{n} r(q) q^{n(n+1)/4}.$    
\end{proposition}

Clearly these two propositions provide a proof of Theorem~\ref{T-precise}.

\begin{proof}[Proof of Proposition~{\rm\ref{P-upper-I}}]
Combining the estimates for $\#\CO(q,n)$ and $\#\CN(q,n)$ given 
in Theorem~1.2 of~\cite{dipippo-howe}, we find that the quantity
$\#\CI(q,n)  - v_{n} r(q) q^{n(n+1)/4} $
is less than or equal to
$$6^{n^2} c_1^n c_2 \frac{n(n+1)}{(n-1)!} q^{n(n-1)/4} + 
\left(v_n+6^{n^2} c_1^n c_3 \frac{n(n+1)}{(n-1)!}\right) q^{(n+2)(n-1)/4},$$
so
$$\frac{\#\CI(q,n)}{v_{n} r(q) q^{n(n+1)/4}}
  \le 1 + \frac{c_{2} G_{n}}{c_{3} r(q) q^{n/2}} + 
       \frac{1}{r(q) q^{1/2}}  +
       \frac{G_{n}}{r(q) q^{1/2}}.$$
An easy induction shows that $G_{n} > 2$, and certainly $c_{2}/c_{3} < 2.5$,
so we have 
\begin{align*}
\frac{\#\CI(q,n)}{v_{n} r(q) q^{n(n+1)/4}} 
  & < 1 + \frac{G_{n}}{r(q) q^{1/2}}
      \left(\frac{c_{2}}{c_{3}} + \frac{1}{2} + 1\right)\\
  & < 1 + \frac{4 G_{n}}{r(q) q^{1/2}}.
\end{align*}
Since $q > M_{n,\eps}$ we have $q^{1/2} > 8G_{n}m_{n,\eps}/\eps$,
and combining this with the fact that $r(q) \ge 1/2$ we find that 
$$\frac{\#\CI(q,n)}{v_{n} r(q) q^{n(n+1)/4}}
  \le 1 + \frac{\eps}{m_{n,\eps}}.$$
But $m_{n,\eps}$ is greater than $8$ for $\eps \le 1$, 
so the right-hand side is at most~$1 + \eps/8$.  This proves the
inequality of the proposition.
\end{proof}

Our proof of Proposition~\ref{P-lower-O} is based upon 
Lemmas~\ref{L-good-Weil-numbers} and~\ref{L-good-polynomials}.
We will compute 
\begin{enumerate}
    \item[$\bullet$] a lower bound on the number of degree-$n$ polynomials 
    satisfying hypotheses (2), (4), and (5) 
    of Lemma~\ref{L-good-polynomials},
    \item[$\bullet$] an upper bound on the number of degree-$n$ polynomials 
    satisfying hypothesis (2) but failing hypothesis (3)
    of Lemma~\ref{L-good-polynomials}, and
    \item[$\bullet$] an upper bound on the number of degree-$n$ polynomials 
    satisfying hypothesis (2) but failing hypothesis (1)
    of Lemma~\ref{L-good-polynomials}.
\end{enumerate}
Subtracting the sum of the latter two estimates from
the first estimate will give us a lower bound on the number 
of degree-$n$ polynomials satisfying all the hypotheses of 
Lemma~\ref{L-good-polynomials}.  By Lemma~\ref{L-good-Weil-numbers},
this lower bound will also be a lower bound on $\#\CO_{\asimp}(q,n)$.
The computation of the lower bound on the number of polynomials satisfying
hypotheses (2), (4), and (5) of Lemma~\ref{L-good-polynomials}
will depend on the following lemma, whose proof we will postpone
until the next section.

\begin{lemma}
    \label{L-reduction}
    Let $n>2$ be a positive integer, let $\eps$ be a real number between $0$
    and~$1$, and let $m = m_{n,\eps}$ be as defined at the beginning of
    this section.  
    Then there are at least $m^{n} (1 - \eps/4)$ monic degree-$n$
    polynomials in $(\BZ/m\BZ)[x]$ such that 
    \begin{enumerate}
	\item[(1)]  there exists a prime divisor $p_{1}$ of $m$ such that the
	reduction of $g$ modulo $p_{1}$ is irreducible, and
	\item[(2)]  there exists a prime divisor $p_{2}$ of $m$ such that the 
	reduction of $g$ modulo $p_{2}$ is a linear times an irreducible.
    \end{enumerate}
\end{lemma}

Before we proceed to the proof of Proposition~\ref{P-lower-O} we should
mention a basic correspondence that we will use repeatedly in
our argument.  Fix our prime power~$q$.
Suppose $g$ is a monic polynomial of degree $n$ with integer coefficients,
say 
$$ g = x^{n} + b_{1}x^{n-1} + \cdots + b_{n},$$
and let $f$ be the polynomial defined by $f(x) = x^{n} g(x+q/x)$,
so that 
$$f = (x^{2n} + q^{n})  + a_{1} (x^{2n-1} + q^{n-1}x) + \cdots 
     + a_{n-1}(x^{n+1} + q x^{n-1}) + a_{n}x^{n}$$
for some integers~$a_{i}$.
Then the linear map $\Omega$ from $\BZ^{n}$ to $\BZ^{n}$ that sends
a vector $\bb = (b_{1},\ldots,b_{n})$ to the vector $\ba = (a_{1},\ldots,a_{n})$
is invertible --- in fact, it is represented by a matrix with integer
entries that is lower-triangular with $1$'s on the diagonal.
Thus, if we let $\bb$ range over a set of vectors that reduces modulo
some integer $m$ to the entire set $(\BZ/m\BZ)^{n}$, then 
$\Omega(\bb)$ will also range over such a set, and conversely,
if $\ba$ ranges over such a set, then so will~$\Omega^{-1}(\ba)$.

Note that if $g$ and $f$ are related as above, then $g$ satisfies hypothesis (2)
of Lemma~\ref{L-good-polynomials} if and only if the roots of $f$ 
in the complex numbers all have magnitude $q^{1/2}$ and the roots of $f$
in the real numbers all have even multiplicity.
Also, the roots of $f$ meet this last condition if and only if the
vector
$(a_{1} q^{-1/2}, a_{2}q^{-1}, \ldots, a_{n}q^{-n/2})$
lies in the region $V_{n}$ of $\BR^{n}$ defined in~\cite{dipippo-howe}.
Thus we will be interested in estimating the sizes of the intersections
of certain lattices with~$V_{n}$.

Let $\be_1$, \dots,~$\be_n$ denote the standard basis vectors of~$\BR^n$.
Our arguments will involve two lattices in $\BR^n$: The first lattice,
denoted~$\Lambda$, is generated by the vectors $q^{-i/2}\be_i$, and the second,
denoted~$\Lambda'$, is generated by the same set of vectors, except with
$q^{-n/2}\be_n$ replaced with~$pq^{-n/2}\be_n$.
Thus~$\Lambda\supset\Lambda'$.

\begin{proof}[Proof of Proposition~{\rm\ref{P-lower-O}}]
Let $m = m_{n,\eps}$ and let $\Lambda''$ denote the lattice~$m\Lambda$.
If $\ell$ is a point in $\Lambda''$ let $B_{\ell}$ denote the ``brick''
$$\ell + \{(x_{1},\ldots,x_{n})\in\BR^{n}
     \mid \forall i: 0 \le x_{i} < m q^{-i/2}\be_{i}\}.$$
Let $S$ denote the set of all $\ell\in\Lambda''$ such that 
$B_{\ell}\subseteq V_{n}$.  The proof of Proposition~2.3.1 of~\cite{dipippo-howe}
(see especially p.~435) shows that 
$$\#S\ge \frac{\text{volume $V_{n}$}}{\text{covolume $\Lambda''$}}
      - 6^{n^{2}} c_{1}^{n} c_{3}\frac{n(n+1)}{(n-1)!} 
            \frac{d}{\text{covolume $\Lambda''$}}$$
where $d$ is the {\it mesh\/} of $\Lambda''$ (see p.~434 of~\cite{dipippo-howe}),
which is~$m q^{-1/2}$.  Since the covolume of $\Lambda''$ is
$m^{n} q^{-n(n+1)/4}$, we find that 
$$\#S \ge m^{-n} v_{n} q^{n(n+1)/4} 
      - 6^{n^{2}} c_{1}^{n} c_{3}\frac{n(n+1)}{(n-1)!} 
            m^{-n+1} q^{(n^{2}+n-2)/4}.$$
Thus 
$$m^{n}\#S \ge v_{n} q^{n(n+1)/4}  (1 - mG_{n}q^{-1/2}),$$
and using the fact that $q>M_{n,\eps}$ we find that 
$$m^{n}\#S \ge v_{n} q^{n(n+1)/4} (1 - \eps/8).$$

Now suppose $\ell$ is a lattice point in~$S$, and consider a typical
element $\bx = (a_{1} q^{-1/2}, a_{2}q^{-1}, \ldots, a_{n}q^{-n/2})$
of~$\Lambda\cap B_{\ell}$.  As $\bx$ ranges over all of $\Lambda\cap B_{\ell}$,
the vector $\ba = (a_{1},\ldots,a_{n})$ ranges over a set of $m^{n}$ elements
of $\BZ^{n}$ that reduces modulo $m$ to all of~$(\BZ/m\BZ)^{n}$.
Lemma~\ref{L-reduction} above shows that of the $m^{n}$ polynomials $g$
we obtain from the vectors $\Omega(\ba)$ arising from elements of
$\Lambda\cap B_{\ell}$, at least $m^{n}(1 - \eps/4)$ satisfy hypotheses
(4) and (5) of Lemma~\ref{L-good-polynomials}. 
So for each element of $S$ we obtain at least $m^{n}(1 - \eps/4)$
polynomials satisfying hypotheses (2), (4), and (5)
of Lemma~\ref{L-good-polynomials}. 
Thus the total number of such polynomials is at least 
$m^{n}\#S(1-\eps/4)$, and by the results of the preceding paragraph
this number is at least
$$v_{n}q^{n(n+1)/4}(1-\eps/4)(1-\eps/8),$$
which is greater than $v_{n}q^{n(n+1)/4}(1-3\eps/8).$

Next we estimate the number of polynomials $g$ that satisfy hypothesis
(2) of Lemma~\ref{L-good-polynomials} but that fail to satisfy 
hypothesis (3).  There is a bijection between the set of such polynomials
and the set $\Lambda'\cap V_{n}$, and Proposition~2.3.1 of~\cite{dipippo-howe}
gives upper and lower bounds for the size of the latter set;
in particular, we find that the number of such polynomials differs from 
$(1/p) v_{n} q^{n(n+1)/4}$ by at most
$$\frac{q^{n(n+1)/4}}{p q^{1/2} } 6^{n^{2}} c_{1}^{n} c_{3}\frac{n(n+1)}{(n-1)!},$$
which is $v_{n} q^{n(n+1)/4}G_{n}/pq^{1/2}$.  Since $q$ is at least $M_{n,\eps}$,
this last quantity is at most $v_{n} q^{n(n+1)/4}\eps / (4p m)$, 
which is less than $v_{n} q^{n(n+1)/4}\eps/32$, because
$m>4$ when $\eps\le 1$.  Thus the number of polynomials that
satisfy hypothesis (2) but not hypothesis (3) is at most
$$v_{n}q^{n(n+1)/4}\left(\frac{1}{p} + \frac{\eps}{32}\right).$$

Finally we estimate the number of polynomials $g$ that satisfy hypothesis
(2) of Lemma~\ref{L-good-polynomials} but that fail to satisfy 
hypothesis (1).  Now, a polynomial $x^{2n} + ax^{n} + q^{n}$
has all of its roots on the circle $|z| = q^{1/2}$ if and only
if $|a|\le 2q^{n/2}$, so there are at most $4q^{n/2} + 1$ polynomials meeting
hypothesis (2) but failing hypothesis (1).  It is very easy to 
show that $4q^{n/2} +1 < v_{n} q^{n(n+1)/4}\eps/32$ when~$q>M_{n,\eps}$.

Now, the number of polynomials meeting all five hypotheses of 
Lemma~\ref{L-good-polynomials} is at least as large as the number
that satisfy hypotheses (2), (4), and (5), less the number that
satisfy hypothesis (2) but that fail either hypothesis (1) or hypothesis~(3).
We find that the number of polynomials meeting all five hypotheses
is at least
\begin{multline*}
  v_{n} q^{n(n+1)/4} \left(1 - \frac{3\eps}{8}\right)
- v_{n} q^{n(n+1)/4} \left(\frac{1}{p} + \frac{\eps}{32}\right)
- v_{n} q^{n(n+1)/4} \frac{\eps}{32} \\
 = v_{n} q^{n(n+1)/4}\left(r(q) - \frac{7\eps}{16}\right)
 \ge v_{n} r(q) q^{n(n+1)/4}\left(1  - \frac{7\eps}{8}\right)
\end{multline*}
and this is the statement of Proposition~\ref{P-lower-O}.
\end{proof}

\section{Proof of Lemma~\ref{L-reduction}}
\label{S-reduction-proof}

In this section we will prove Lemma~\ref{L-reduction}.  We continue to
use the notation set at the beginning of Section~\ref{S-asymptotic-higher}.

For the moment, 
let us write $A_{n,p}$ for the set of monic degree-$n$ 
irreducible polynomials in $\BF_{p}[x]$
and $B_{n,p}$ for the set of monic degree-$n$ polynomials in $\BF_{p}[x]$
that factor as a linear polynomial times an irreducible.
\begin{lemma}
    \label{L-polynomial-count}
    Let $p$ be a prime.  For all $n>0$ we have
    $\#A_{n,p} \ge p^{n}/(2n)$,
    and for all $n>1$ we have 
    $\#B_{n,p} \ge p^{n}/(2n-2)$.
\end{lemma}

\begin{proof}
The lemma follows easily from the well-known exact formula
$$\#A_{n,p} = \frac{1}{n} \sum_{d\mid n} p^{d}\mu\left(\frac{n}{d}\right),$$
where $\mu$ is the M\"obius function.
\end{proof}

Suppose~$n > 1$.
We see from Lemma~\ref{L-polynomial-count} that if we choose a monic degree-$n$
polynomial $f$ at random from $\BF_{p}[x]$ (with the uniform distribution),
then 
\begin{align*}
\Prob(f\not\in A_{n,p}) & \le 1 - \frac{1}{2n} \\
\intertext{and}
\Prob(f\not\in B_{n,p}) & \le 1 - \frac{1}{2n-2}.
\end{align*}

Suppose that $\eps$ is given, with~$0 < \eps < 1$.
Let $k = k_{n,\eps}$ and $m = m_{n,\eps}$ be as at the beginning of 
Section~\ref{S-asymptotic-higher}, so that $m$ is the product of
the first $k$ prime numbers.  Now suppose we choose a monic degree-$n$
polynomial $f$ at random from~$(\BZ/m\BZ)[x]$.  By the Chinese remainder
theorem, making such a choice is equivalent to choosing a monic degree-$n$
polynomial $f$ at random from $\BF_{p}[x]$ for each of the first $k$
primes~$p$.  Thus we see that
\begin{align*}
\Prob(\forall p\mid m : (f\bmod p) \not\in A_{n,p})
     & \le \left(1 - \frac{1}{2n}\right)^{k} \\
\Prob(\forall p\mid m : (f\bmod p) \not\in B_{n,p})
     & \le \left(1 - \frac{1}{2n-2}\right)^{k},
\end{align*}
and it follows that
\begin{multline*}
\Prob(\exists p_{1}, p_{2} \mid m : (f \bmod p_{1}) \in A_{n,p_{1}}
                          \text{\ and\ } (f \bmod p_{2}) \in B_{n,p_{2}}) \\
	> 1 - \left(1 - \frac{1}{2n}\right)^{k}
	    - \left(1 - \frac{1}{2n-2}\right)^{k} 
	> 1 - 2\left(1 - \frac{1}{2n}\right)^{k}.
\end{multline*}
But the definition of $k_{n,\eps}$ shows that 
$$\left(1-\frac{1}{2n}\right)^{k}  \le \frac{\eps}{8},$$
so 
$$
\Prob(\exists p_{1}, p_{2} \mid m : (f \bmod p_{1}) \in A_{n,p_{1}}
                          \text{\ and\ } (f \bmod p_{2}) \in B_{n,p_{2}}) 
	> 1 - \eps/4.
$$
Thus the number of monic degree-$n$ polynomials in $(\BZ/m\BZ)[x]$
that satisfy the two conditions of  Lemma~\ref{L-reduction}
is at least $m^{n}(1-\eps/4)$, as was to be shown.

\end{document}